\documentclass[a4paper,10pt]{article}
\usepackage{graphicx}
\usepackage{epsfig}
\usepackage{fancyhdr}

\usepackage{amsmath,amsfonts,amssymb,amsthm}
\usepackage{setspace}
\theoremstyle{plain}

\newtheorem*{thm*}{Theorem}

\newtheorem*{lem*}{Lemma}

\theoremstyle{definition}

\newcommand{\reftit}{\textit}    
\newcommand{\refis}{\textbf}     

\begin{document}

\title{Mixing times via super-fast coupling}
\author{Robert M. Burton \footnote{
 Department of Mathematics,
 Oregon State University,
 Corvallis, OR
 97331-4605, USA
 \texttt{bob@oregonstate.edu} }
 \and Yevgeniy Kovchegov \footnote{
 Department of Mathematics,
 Oregon State University,
 Corvallis, OR
 97331-4605, USA
 \texttt{kovchegy@math.oregonstate.edu}}}
\date{ }
\maketitle

 \begin{abstract}
 For the probabilistic model of shuffling by random transpositions we provide a coupling construction with the expected coupling time of order $C n \log (n)$, where $C$ is a moderate constant. We enlarge the methodology of coupling by including intuitive non-Markovian  coupling rules. We discuss why a typical Markovian coupling is not always sufficient for resolving mixing time questions.
 \end{abstract}


\section{Introduction}
Random shufflings of a deck of $n$ distinct cards are well-studied objects, and a frequent metaphor describing a class of Markov chains invariant with respect to the symmetric group, $S_n$. Here the focus is on transposition shuffling, one of the simplest shuffles, defined by uniformly sampling the deck twice with replacement, and then interchanging the positions of these cards, if they are different.

Clearly, as it is always the case for mixing finite state Markov chains, the distribution of the ordering of the deck converges in total variation to the invariant measure, which by symmetry is uniform on the permutation group, $S_n.$  It is the rate of mixing which presently holds our interest, as well as the coupling methods by which one might attempt to show good upper bounds on this rate.

This is a well-defined problem in probability theory, and one would expect that a coupling argument would be the instrument of first choice. Indeed, there is such an approach, given in the online notes of Aldous and Fill \cite{af} (also see \cite{a}). This method gives a rate of $O(n^2)$, and will be discussed later in the article.  Unfortunately, but necessarily, this rate is not the optimal rate [we  expect, $O(n\log{n}),$ which was proven by Diaconis and Shahshahani \cite{ds} using methods from 
representation theory]. This gap is apparent, and somewhat long-standing;
indeed, Peres has listed the problem of showing the $O(n\log{n})$ as the rate of uniform mixing, using a coupling approach, as one of a number of interesting open problems (see \cite{p}). The problem is also mentioned in Saloff-Coste \cite{sc} and other important publications. Here we solve this problem. Moreover, we deconstruct different ways of looking at this kind of problem, and enlarge the intuitions that might guide us in coupling.

We would like to mention that Matthews \cite{m} gives a purely probabilistic proof of the $O(n \log{n})$ bound using strong stationary
times.  Diaconis and Saloff-Coste \cite{dsc} also treat a host of other transposition problems, and a survey of the many appearances
of random transpositions can be found in Diaconis \cite{d03}. Among the publications related to this research, we would like to list the recent work of Berestycki and Durrett \cite{bd} on random transposition, and the paper of Hayes and Vigoda on a non-Markovian coupling for graph coloring problems.

 Coupling comes from at least two fields. Of course, one is probability theory where, typically, the intent of the proof is to give an intuition of the phenomena and why it happens. For this reason, most proofs are Markovian and depend solely on information available to the present moment of the evolution of a Markov process or process derived from a Markov process such as a weak Bernoulli process. Then one could, for example, use coupling to show strong versions of the central limit theorem (see for example, \cite{bbd}). Another field that coupling developed within is ergodic theory. Here, the problems had complicated forms of dependence and the coupling was non-intuitive, using the marriage lemma to their existence. See \cite{o} and \cite{fo}. The result that any two weakly Bernoulli processes with the same entropy are isomorphic in a measure-theoretic translation invariant way was proven in such a way. Indeed, the choice of coupling depends essentially on the distinction between quantitative and qualitative. To show the ideas one would prefer Markovian coupling, even though this kind of coupling cannot be best possible and to get precise rates one would use non-Markovian coupling because no other kind would work.
 The history of research on non-Markovian coupling techniques and maximal coupling of Markov chains goes back to the publications of D.Griffeath, J.Pitman, and S.Goldstein in the 70's (see \cite{dg}, \cite{jp}, and \cite{sg}). K.Burdzy and W.S.Kendall have studied the efficient Markovian couplings in \cite{kbwk}.

\subsection{Preliminary details}
 Suppose that the number of cards $n$ is fixed and define ${\cal S} = \{0, 1, \dots, n-1\}$, and suppose ${\cal S}^2$ is the state space of a sequence of i.i.d. random pairs  $<{\bf a},{\bf b}>$ which are uniformly distributed on ${\cal S}^2$. Each  pair $<{\bf a},{\bf b}>$ generates an independent transposition, although there is a probability of $\frac{1}{n}$ that there will be no change in the ordering of the deck. Given an initial probability distribution on $S_n$, we let the law of the distribution on $S_n$ be actualized by a random permutation $X_t$ describing the configurations of cards. This gives a well defined Markov chain. It will sometimes be useful to Poissonify the Markov chain to continuous time with i.i.d. exponentially distributed
 interarrival times. 
 
 \vskip 0.2 in
 \noindent
A coupling argument would require the construction of a joint distribution for $X_t$ and $Y_t$ on the product of sample spaces $S_n \times S_n$, where the distribution of $X_0$, $dist(X_0)$ is the point mass distribution on $S_n$, giving the identity configuration probability 1. This joint distribution must have correct marginal distributions, and the random permutations $X_t$ and $Y_t$ must agree from some time $T$ onward.
Notice that the shuffle is invariant with respect to the symmetric group in the following sense:
\vskip 0.2 in
{\bf Definition.}   Suppose that ${X_t}$ is a Markov Chain on $S_n$. The chain is said to be {\bf group invariant} if  $dist(\gamma X_{t+1} | X_t = \alpha)$ is equal to the $dist( X_{t+1} | X_t = {\gamma}^{-1}{\alpha})$ for all $\gamma, \alpha \in S_n$.
\vskip 0.2in

This is a homogeneity condition similar to that of independent increments, and says that the shuffle is independent of the values or printed labels on the cards. It implies that the distribution of cycle structures of $X_1$ given $X_{0} = \alpha$ depends only on the cycle structure of $\alpha$. The set of group elements with a given cycle structure form a conjugacy class in the group, $S_n$. For example, if $\alpha$ is a transposition, i.e. a 2-cycle, then regardless of what $\alpha$ is specifically it has the same coupling and transition structure as any other transposition. Since the identity transitions to each of the $\begin{pmatrix}  n\\ 2 \end{pmatrix}$ transpositions with equal probability $\frac{2}{n^2}$, and remains at the identity with probability $\frac{1}{n},$ we could analyze this process as a Markov chain on the conjugacy classes, i.e. classes of permutations with the same cycle structure.  Further, this random walk on the cyclic decompositions is biased toward intermediate central weights of the permutations (where the weight is equal to the sum of the cycle lengths minus the number of cycles, and equal to the minimum number of transpositions required to reduce the permutation to the identity). We could have used this approach to analyze recurrence rates. The number of distinct cycle structures is equal to the partition function on the integers which grows exponentially with rate $\sqrt{\frac{4n}{3}}$ much slower than the rate of $n\log(n)$ that $n!$ grows with.

There are many ways to associate an element $<{\bf a},{\bf b}>$ in ${\cal S}^2$ with a transposition. Each card has a value or label printed on it, and each card has a location: the number of cards above it on the deck plus 1. Let $Q$ be the set of values or labels, and let $P$ be the set of positions of the cards. If $<{\bf a},{\bf b}>$ is a member of $P\times P$ then we associate $<{\bf a},{\bf b}>$ with interchanging the locations of a card at  location ${\bf a}$ and a card at location ${\bf b}$. If $<{\bf a},{\bf b}>$ is a member of $Q \times P$ (in which case we will use $<\fbox{\bf a},{\bf b}>$ notation), then we associate $<\fbox{\bf a},{\bf b}>$ with taking card $\fbox{\bf a}$ and placing it at location ${\bf b}$, and then using the card that was in location ${\bf b}$  to replace $\fbox{\bf a}$ in its original location. When we begin to couple, we will be following the motion of two decks. It is natural to apply the same operation to each one of the two decks. For example, if we did this with $Q\times P$ association for evolution, we could do the same transposition $<\fbox{\bf a},{\bf b}>$ for both to eventually obtain coupling in $O(n^2)$, as in Aldous and Fill \cite{af}. We however have much more flexibility than that with a bijective map from $Q \times P$ to $Q \times P$, called an {\bf association mapping} to relate the coupled moves of each deck. An association mapping tells us how to couple the immediate descendants of two group elements that are coupled. Recall that we must preserve relations, like siblings and cousins in the tree.

Most coupling arguments are made up of present and past measurable constituents; that is the coupling method is adapted to an increasing sequence of $\sigma$-fields so that the present and past are measurable with respect to their corresponding $\sigma$-fields, and that the $\sigma$-fields have enough extra randomness to perform independent experiments, subdivide atoms, and so on. In some ways this tendency toward adaptive coupling is historical, and in some ways it is natural to follow one's intuition, and then make it rigorous. The first coupling arguments most people see is a passive coupling in which two Markov chains are allowed to go their own ways independently of each other until they happen to obtain the same state at the same time. From this random time onward they are coupled together. It is also true that if we had perfect knowledge of the situation, we would be able to make the optimal coupling at any time. Unfortunately, the numbers are usually too large and the relationships too complex. The point here is that non-adapted coupling can be natural, intuitive, and have great power in solving problems. Indeed, as we shall see, we can even couple so that we anticipate the future and prepare for it, while maintaining the only essential ingredient of coupling, that of having perfect (or near perfect) distributions on the marginal processes. Moreover, it is possible to have couplings made up of surgically cut and pasted pieces of sample path.

\subsection{Strong uniform mixing, weak Bernoulli, and coupling}

  If ${X_t}$ and ${Y_t}$ are stochastic processes, then a coupling of ${X_t}$ and ${Y_t}$ is a joint probability distribution on the product, ${(X_t, Y_t)}$ so that the marginal distributions on $X_t$ and $Y_t$ agree with original distributions. We say that $X_t(\omega)$ and $Y_t(\omega)$  are coupled at the random time $T$ if for $t \ge T$ it is the case that $X_t = Y_t$. If $T$ is finite a.s., then this argument shows that the distributions of $X_t$ and $Y_t$ are converging in the total variation distance as $t$ becomes large.

\vskip 0.2 in
\noindent
  Diaconis and Shahshahani \cite{ds} define a finite state Markov chain ${X_t}$ with invariant probability $U$ to be {\bf strong uniform mixing} if there is a stopping time $T$, so that $P[ T = t, X_t = \alpha]$ is independent of the group element $\alpha$. Because of this and the invariance of the uniform distribution on the group $G$, for any $\alpha$ we have $P[ X_t = \alpha | t \ge T] = 1/|G|$, as the invariant measure is uniform on $G$. Regardless of the distribution of  $X_0$, we have the following bound on the total variation norm, $||U - dist(X_t)||_{TV} \leq P[T > t]$, because this is the only part of the probability space where the total variation is not forced to be 0. The rate of mixing is carried by the distribution of the stopping time.
Coupling arises because the total variation norm is achieved by joining the distributions together in a probability preserving way. The process is Markov, so it makes sense for the definition to be independent of the initial state. Coupling for general processes is usually connected to weak Bernoulli, which has a rich history (a.k.a. absolutely regular and $\beta-$mixing \cite{vr}).
\vskip 0.2in
{\bf Definition.} A finite-valued stochastic process ${X_t}$ is {\bf weak Bernoulli} if there is a coupling  $\{(X'_t,X''_t) : t \in \mathbb{Z}\}$ such that
(i) $\{X'_t\}$ and $\{X''_t\}$ have the same distribution as $\{X_t\}$, (ii) the past of $X''_t$, $\{X''_t : t = 0,-1,-2, \ldots \}$ is independent of $\{X'_t\}$, and (iii) there is a random a.s. finite time $T$ so that $t \ge T$ implies that $X'_t = X''_t$.
\vskip 0.2in

In this case, the future becomes independent of past values in a strong way.
The pathwise coupling version of weak Bernoulli is

\vskip 0.2in
{\bf Definition.} A finitely valued stochastic process ${X_t}$ is {\bf tree weak Bernoulli} if there is a coupling  $\{(X'_t,X''_t) : t \in \mathbb{Z}\}$ such that
(i) $\{X'_t\}$ and $\{X''_t\}$ have the same distribution as $\{X_t\}$, (ii) the past of $X''_t$, $\{X''_t : t = 0,-1,-2, \ldots \}$ is independent of $\{X'_t\}$,  (iii) there is a random a.s. finite time $T$ so that $t \ge T$ implies that $X'_t = X''_t$, and (iv) the coupling respects the tree structure of future sample paths, so if $X'_t$ and $X''_t$ are coupled at time $t>0$ then each of their descendants (or successors) are also coupled together, i.e. the coupling is given by a tree automorphism of the branching future paths.
\vskip 0.2in
\noindent
The terminology comes from Hoffman and Rudolph \cite{hr}, in which they used tree very weak Bernoulli to study isomorphisms of $1$ to $p$ endomorphisms. The conditional distributions of $\{X_t\}_{t \in \{1,2, \ldots, M\} }$ live on the set of labeled trees of length $M$. We raise the definition because many, if not most, coupling methods have this property of being tree consistent.

When applying tree weak Bernoulli or tree coupling, it is important to consider the events defined by the finitely valued random
variables to be atomless, e.g. a subset of the unit interval. Achieving optimal total variation norm typically requires subdividing
these events into events of smaller but arbitrary probability.

In any given situation, it is possible, in principle, to find the optimal coupling because the total variation norm is achievable. However, the size of the state space and dependence within the process often make this impractical.

We developed methods to speed up the coupling time while maintaining the distributions on the marginal
processes. The methods rely on intuition and insight, and were usually dynamical. By dynamical, we mean that
the coupling at a given time depended upon the process up to that time, and, perhaps, some external
sources of independent randomness. In other words, the coupling was adapted to the pair of processes as they were
being constructed. Unlike uniform strong mixing, tree weak Bernoulli is a path-wise  property and requires a probability-preserving
path-wise isomorphism between the future trees of possibilities of the processes.

To illustrate the difference between uniform strong mixing and tree weak Bernoulli
suppose we have a Markov chain with state space
$\{ a_1, a_2, b_1, b_2 \}$, and suppose the only allowable transitions, each with
conditional probability one-half, are as given in the following transition probability matrix
$$\begin{array}{|c||cccc|} \hline & a_1 & a_2 & b_1 & b_2 \\\hline \hline a_1 & 0 & 0 & 1/2 & 1/2 \\ a_2 & 1/2 & 0 & 1/2 & 0 \\ b_1 & 1/2 & 1/2 & 0 & 0 \\ b_2 & 0 & 1/2 & 0 & 1/2 \\ \hline \end{array} $$

The matrix is doubly stochastic, so the invariant probability
is uniform on each state. Suppose we have two
Markov chains with this law, and initial states
$X_0 =  a_1$, and $Y_ 0 = b_2$. The possible forward paths are

$$\left.\begin{array}{ccccccccc}  &   &   &   & a_1 &   &   &  &   \\   &   &   & \swarrow &   & \searrow &   &   &   \\
&   & a_2 &   &   &   & b_1 &   &   \\  &  \swarrow & \downarrow &   &   &   & \downarrow & \searrow   &  \\  b_1  &   & b_2 &   &   &   & a_1 &   &  a_2 \\\end{array}\right.
\quad \text{ and } \quad
\left.\begin{array}{ccccccccc}  &   &   &   & b_2 &   &   &  &   \\   &   &   & \swarrow &   & \searrow &   &   &   \\
&   & a_1 &   &   &   & b_2 &   &   \\  &  \swarrow & \downarrow &   &   &   & \downarrow & \searrow   &  \\  a_2  &   & b_1 &   &   &   & a_1 &   &  b_2 \\
 \end{array}\right.$$
 The distribution of $X_2$ and $Y_2$ are equal and at
equilibrium, $|| dist(X_2) - dist(Y_2)||  =
||dist(X_2) - U|| = 0$ while $|| dist(X_1) - dist(Y_1)||  =
1$ so $X_1 \ne Y_1$ a.s. On the other hand, any tree coupling
of $(X_0,  X_1, X_2)$ with $(Y_0,  Y_1, Y_2)$ has
$P[ X_2  =  Y_2] = \frac{1}{2}.$

Similarly, the tree weak Bernoulli, or tree coupling coefficients of a
process, may be different than the uniformly strong mixing
coefficients, although the coupling distance is never less
than the total variation distance.

\section{Transposition shuffling for the case $n =3$}

 The models we deal with here are invariant random walks on a group in which each step has the same probability, although the case that different steps will go  to the same location is allowed and expected. The set up is a Markov chain that is itself a hidden Markov chain with the same number of equally likely outcomes at any stage.

  We begin with the group $S_3$. We will make some comparisons between two processes, $(X_t)$ where $X_0 = 123$,  and $Y_t$.
  Consider the adjacency matrix, $M$,
 $$\begin{tabular}{|c||cccccc|}\hline  & 123 & 231 & 312 & 132 & 213 & 321 \\\hline\hline 123 & 3 & 0 & 0 & 2 & 2 & 2 \\ 231 & 0 & 3 & 0 & 2 & 2 & 2 \\ 312 & 0 & 0 & 3 & 2 & 2 & 2 \\ 132 & 2 & 2 & 2 & 3 & 0 & 0 \\ 213 & 2 & 2 & 2 & 0 & 3 & 0 \\ 321 & 2 & 2 & 2 & 0 & 0 & 3 \\\hline \end{tabular}$$
 If $M$ is the matrix shown then the $(1/9)M$ is the stochastic matrix that represents the dynamics of the chain.

 The square of $M$ is $M^2 = 9I + 12N$ where $N$ is the 6 by 6 matrix consisting of all 1s. The maximum total variation distance between any starting point and any other starting point is $(21 - 12)/81 = 1/9$. If we raise $M$ to the power $2t$ we get  $M^{2t} = (9I + 12N)^t = a_tI + b_tN$, where $a_t=9^t $. This means that after $2t$ iterations the total variation norm between any two distinct starting places is $1/3^{2t}$,  so the mixing rate is geometric with error of ${1 \over 2 \cdot 3^{2t}}$. The 1/2 comes from reducing the total variation distance between distinct starting places to that of one starting permutation and the uniform distribution.

Since every odd power of $M$ is the product of an even power of $M$ with $M$, we can similarly compute the total variation distance in this case. An analogous calculation gives
the total variation distance between
$X_t$ and the uniform $U$ as
$$||dist(X_t) - U||_{TV}  = \begin{cases}
      {5 \over 2} \cdot 3^{-t} & t \text{ odd }, \\
      {1 \over 2} \cdot 3^{-t} & t \text{ even }
\end{cases}$$
So for all $t$ we have
$${1 \over 2} \cdot 3^{-t} \le  || dist (X_t) -U||_{TV} \le {5 \over 2} \cdot 3^{-t}$$

If we use the $Q\times P$ association, picking a face value and then picking a location to swap cards with, and if $X_0 = 123$, then the first transposition will have nine equally likely outcomes $X_1$, three of which are configuration 123, and there are two chances for each of three configurations, 132, 213, and 321. These are shown in the chart below. In keeping with group invariance, each turn there are three chances to remain in its previous state and two chances for each of two group elements of opposite parity.

$$\left.\begin{array}{|c||ccc|}\hline Q \setminus P & {\bf 1} & {\bf 2} & {\bf 3} \\\hline \hline \fbox{\bf 1} & 123 & 213 & 321 \\ \fbox{\bf 2} & 213 &  123 & 132 \\ \fbox{\bf 3} & 321 & 132 &   123 \\\hline \end{array}\right.$$

On the other hand, if we wish to couple with the path arising from $X_0=132$ using group invariance, we get the following table:

$$\left.\begin{array}{|c||ccc|}\hline Q \setminus P & {\bf 1} & {\bf 2} & {\bf 3} \\\hline \hline \fbox{\bf 1} & 132 & 312 & 231 \\ \fbox{\bf 2} & 231 &  123 & 132 \\ \fbox{\bf 3} & 312 & 132 &  123 \\\hline \end{array}\right.$$

Looking at the diagram, we see that four out of nine outcomes can be perfectly coupled, and that the remaining line up as transpositions of each other. Actually, we can do a bit better by coupling two that differ by a 3-cycle and then four that couple as transpositions. Specifically, couple two 123s with the corresponding 123s, the two 132s with the corresponding 132s, the remaining 123 with 312, and 213 with 132, and then the remaining three couple in any way with the corresponding remaining 3 outcomes. After drawing some infinite trees, and calculating many geometric series one arrives at the optimal tree coupling which is $m 3^{-m}$. Note that this coupling is strictly larger, leaving open the possibility that tree coupling may be unable to achieve the optimal mixing rates. As an aside, note that the trajectories of the process are given by a two dimensional substitution system given by iterating the three by three diagrams above.

\section{Coupling to the future: an example}
In this short section we present an example of a coupling that utilizes the same approach as the coupling construction for the shuffling by random transpositions that will be described in section \ref{superfast}.
Consider a continuous-time process $(X_t,Y_t)$ on $S=\{0,1\}\times\{0,1,2\}$ with generator
$$\begin{tabular}{|c||cccccc|}  \hline & (0,0) & (0,1) & (0,2) & (1,0) & (1,1) & (1,2) \\\hline\hline (0,0) & -12 & 0 & 0 & 10 & 1 & 1 \\ (0,1) & 0 & -12 & 0 & 1 & 10 & 1 \\ (0,2) & 0 & 0 & -12 & 1 & 1 & 10 \\ (1,0) & 10 & 1 & 1 & -12 & 0 & 0 \\ (1,1) & 1 & 10 & 1 & 0 & -12 & 0 \\ (1,2) & 1 & 1 & 10 & 0 & 0 & -12 \\\hline \end{tabular}$$
\noindent
  Here the first coordinate $X_t$ changes to $1-X_t$ with rate 10, independently of the second coordinate $Y_t$. The second coordinate switches
  $~~~Y_t \rightarrow (Y_t-1) \!\!\mod 3~~~$ or  $~~~Y_t \rightarrow (Y_t+1) \!\!\mod 3~~~$ with rate 1 each. However every time the second coordinate changes, the first coordinate must also change.
\vskip 0.2 in
\noindent
 We want to find a fast coupling for  the above process. Coupling the first coordinate is simple:  wait with rate 20, then assign with
 equal probabilities of ${1 \over 2}$ either 0 or 1
 to both $X_t$ and $X'_t$. Coupling the second coordinate is similarly simple: wait with rate 3 before assigning any one of the three
 values 0, 1, or 2 to both $Y_t$ and $Y'_t$. However coupling both coordinates simultaneously for the processes $(X_t,Y_t)$ and $(X'_t,Y'_t)$
 may create the following complication: let $T_1$ be the exponential r.v. with parameter 20, and $T_2$ be the exponential r.v. with parameter 3.
 With large probability, the first coordinates will couple before the second, i.e. $T_1<T_2$. In two out of three cases, either $Y_t$ or $Y'_t$
 will not change at $T_2$. Therefore in two thirds of the cases, when we couple the second coordinate, we simultaneously decouple the first. Thus
 we will need extra time for the first coordinates to couple again.
\vskip 0.2 in
\noindent
  Coupling to the future works in the following way. We start by generating $T_1$ and $T_2$. If $T_1<T_2$, we will start by deciding the value
  for $Y_{T_2}=Y'_{T_2}$. If it is different than both $Y_0$ and $Y'_0$, then we need to randomly generate $X_{T_1}=X'_{T_1}$, and
  we are done, the process is coupled at $T_2$. However, if $T_1<T_2$ and $Y_{T_2}=Y'_{T_2}$ matches either $Y_0$ or $Y'_0$,
  we will need to pair $X_{t}$ and $X'_{t}$ differently. Namely we will need to couple $X_{T_1}$  and $1-X'_{T_1}$, and keep
  $X_t=1-X'_t$ paired until time $T_2$. Therefore at $T_2$, when the second coordinates $Y_t$ and $Y'_t$ couple, one of the first coordinates
  must also flip, thus coupling the two processes, $(X_t,Y_t)$ and $(X'_t,Y'_t)$. In any situation, the coupling time will be $\max\{T_1,T_2\}$, while
  before it was larger.
\vskip 0.2 in
\noindent
  Looking into the future of $Y_t$ and $Y'_t$, we decided whether to pair $X_t$ with $X'_t$ or $X_t$ with $1-X'_t$.
  We will call this pairing an {\it association map}.

\section{Super-fast coupling} \label{superfast}
 We consider an $n$ card deck. A random transposition shuffle is generated by making two independent uniform choices of
 cards, and interchanging them. We assume this being a continuous time process, where the transpositions happen one at a time with exponential waiting times of rate one in between. For each pair of distinct cards $\fbox{\bf a}$ and $\fbox{\bf b}$
 there are two ways to order them, $< \fbox{\bf a}, \fbox{\bf b} >$ and $< \fbox{\bf b}, \fbox{\bf a} >$, each with equal probability.
 Thus, each transposition  $< \fbox{\bf a}, \fbox{\bf b} >$ occurs with rate ${2 \over n^2}$, i.e. the two
 identical transpositions $< \fbox{\bf a}, \fbox{\bf b} >$ and $< \fbox{\bf b}, \fbox{\bf a} >$ will each happen with rate ${1 \over n^2}$. The choice of two identical cards  happens with probability ${1 \over n}$, and the unnecessary transposition $<\fbox{\bf a}, \fbox{\bf a}>$ is therefore assigned the rate ${1 \over n^2}$.
 \vskip 0.2 in
\noindent
 Diaconis and Shahshahani used group representation methods to show that the mixing time
 for this shuffling process is $O(n\log{n})$. See \cite{ds} and \cite{d}.  In \cite{af} a classical coupling approach was shown
 to produce the upper bound of order $O(n^2)$,  while \cite{p} lists showing $O(n\log{n})$ mixing time via {\it coupling} construction
 as an open problem. In this section we produce such coupling construction.
 \vskip 0.3 in
 \noindent
 Throughout this section we will use the following notations and vocabulary.
\begin{center}
\begin{tabular}{|c|p{4.5 in}|}
\hline
\multicolumn{2}{|c|}{\bf Notations and vocabulary} \\
\hline
$<\cdot,\cdot>$ & - transpositions in the card shuffling process\\
$<\fbox{\bf a},\cdot>$ & - transpositions initiated by card $\fbox{\bf a}$\\
$A_t$ & - the top shuffling process\\
$B_t$ & - the bottom shuffling process\\
$\begin{pmatrix}  A_t\\ B_t \end{pmatrix}$ & - the coupled process\\
$<\cdot,\cdot>_A$ & - transpositions in the top shuffling process $A_t$\\
 $<\cdot,\cdot>_B$ & - transpositions in the bottom shuffling process $B_t$\\
 $\ll \cdot, \cdot \gg$ & - simultaneous transpositions in the coupled process $\begin{pmatrix}  A_t\\ B_t \end{pmatrix}$\\
{\it association map} & - association between positions/locations in the top process and positions/locations in the bottom process that will be used to establish the rates for the coupled process\\
\hline
\end{tabular}
\end{center}

\subsection{Label-to-location and  label-to-label transpositions.}
 One of the possible coupling constructions was described in \cite{af}. There, at each step, a card $\fbox{\bf a}$ and a location $i$
 were selected at random, and the transposition $\ll \fbox{\bf a},i \gg$ that moves card $\fbox{\bf a}$ to location $i$ in both
 the top and the bottom processes, $A_t$ and $B_t$, was applied at each iteration. That is, we used $Q\times P$ association for both decks.
 Clearly, this coupling slows down significantly when the number of discrepancies is small enough,
 thus producing an upper bound of order $O(n^2)$, instead of $O(n\log{n})$. 
  \vskip 0.2 in
\noindent
An equivalent upper bound can be achieved by slightly modifying the coupling rules. On each iteration, we
 randomly select a card $\fbox{\bf a}$. If $\fbox{\bf a}$ is not coupled, then we apply transposition $\ll \fbox{\bf a},i \gg$ for a randomly selected location $i$. If $\fbox{\bf a}$ is one of the coupled cards, we apply transposition $\ll \fbox{\bf a}, \fbox{\bf b} \gg$, transposing $\fbox{\bf a}$ with a randomly selected card $\fbox{\bf b}$, in both processes, $A_t$ and $B_t$.
 \vskip 0.2 in
\noindent
 This latter, slightly modified coupling is important, as it can be improved
 to match the correct  $O(n\log{n})$ order for mixing time. The improvement comes in the form of a combinatorial trick similar in spirit to the one
 used by Euler in computing the number of permutations of $n$ elements with all elements displaced.
 \vskip 0.2 in
\noindent
 For the rest of the section, transpositions $\ll \fbox{\bf a},i \gg$ will be called {\it label-to-location}, while transpositions $\ll \fbox{\bf a}, \fbox{\bf b} \gg$ will be called {\it label-to-label}.
  \vskip 0.2 in
\noindent
 Now, without loss of generality, we would like to state that,  in the last coupling construction, whenever a coupled card, for example $\fbox{\bf a}$, is selected, and thus a transposition  $\ll \fbox{\bf a}, \fbox{\bf b} \gg$ with another card $\fbox{\bf b}$ is to be applied, we actually need not
 do this transposition, as the resulting discrepancies will be the same (with different card values, of course) as before the transposition.
 In general, for any permutation $\sigma$ of face values on the cards, the coupled process
 $\begin{pmatrix}  A_t\\ B_t \end{pmatrix}$  is isomorphic to $\begin{pmatrix}  \sigma{A}_t\\ \sigma{B}_t \end{pmatrix}$
 under what we will call {\it the group invariance} of the coupled process.
 The point of mentioning the group invariance property is to stress our right to suppress label-to-label transpositions in the coupled process whenever necessary.

\subsection{Improving the coupling}

The basic strategy for the type of coupling constructions described in this paper is to condition on a key $\sigma$-field in the future and then to use this future information to arrange the intermittent events so that the process is set up for a successful coupling event. This is completely legitimate as long as we take care to have the marginal stochastic processes maintain the correct finite-dimensional   distributions. If this is followed, we are able to modify the joint distributions of the processes however we like to obtain our goal.
\vskip 0.2 in
\noindent
In order to illustrate the forthcoming coupling construction and the notations used therein let us consider the following example.
\vskip 0.2 in
\noindent
{\bf Example.}
 Consider a deck of four cards that are paired so they have two discrepancies ($d=2$) at time $t_0=0$.
   $$\begin{matrix}
   \quad \text{deck }A: \quad \fbox{1} & \fbox{2} & \fbox{3} & \fbox{4} \\
   \quad \text{deck }B: \quad \fbox{1} & \fbox{3} & \fbox{2} & \fbox{4} \\
   \text{location}:\quad ~1\! & 2 & 3 & 4
 \end{matrix}$$
 We pick a random uniform location $i_1 \in \{1,2,3,4\}$ and an exponential time $t_1$.
 \vskip 0.3 in
Conditioned on the event of card $\fbox{2}$ jumping to $i_1$ (i.e. $<<\fbox{2},i_1>>$) at time $t_1$, we provide the following coupling rules for time $t \in [0,t_1]$.
 \vskip 0.2 in
 \noindent
 {\bf Case I: $i_1=2$ or $3$.} \\
 Here cards $\fbox{1}$, $\fbox{3}$ and $\fbox{4}$ do label-to-label jumps only, and each of these jumps leads to card setups that are isomorphic to the original setup, up to relabeling the cards. Because of this we suppress noting any change at all.
 To illustrate this, a label-to-label jump $<<\fbox{1},\fbox{3}>>$will take
   $$\begin{matrix}
    \quad \text{deck }A: \quad \fbox{1} & \fbox{2} & \fbox{3} & \fbox{4} \\
   \quad  \text{deck }B: \quad \fbox{1} & \fbox{3} & \fbox{2} & \fbox{4} \\
 \end{matrix} ~~~\text{ to }~~~
  \begin{matrix}
    \quad \text{deck }A: \quad \fbox{3} & \fbox{2} & \fbox{1} & \fbox{4} \\
   \quad  \text{deck }B: \quad \fbox{3} & \fbox{1} & \fbox{2} & \fbox{4} \\
 \end{matrix}$$
 This latter set up is equivalent to the former setup, so this case leads to no change.
 At the end, card $\fbox{2}$ does the label-to-location jump $<<\fbox{2},i_1>>$ at time $t_1$, canceling the discrepancies.
 \vskip 0.3 in
 \noindent
 {\bf Case II: $i_1=1$ or $4$,} one of the non-discrepancy locations. Again it suffices to only consider $i_1=4$. Here are the rates for
 this case:
\begin{itemize}
\item Cards $\fbox{1}$ and $\fbox{3}$ do label-to-label jumps, and as in Case I, no real change occurs and we again suppress any notational changes.
\item On the other hand if we pick card $\fbox{4}$, then we couple both decks together as follows:
either we interchange cards $\fbox{4}$ and $\fbox{3}$ on the top ($<\fbox{4},\fbox{3}>_A$), getting
   $$\begin{matrix}
    \quad \text{deck }A: \quad \fbox{1} & \fbox{2} & \fbox{4} & \fbox{3} \\
   \quad  \text{deck }B: \quad \fbox{1} & \fbox{3} & \fbox{2} & \fbox{4} \\
 \end{matrix}$$
 The other possibility is we interchange $\fbox{4}$ and $\fbox{3}$ on the bottom ($<\fbox{4},\fbox{3}>_B$) to get
    $$\begin{matrix}
    \quad \text{deck }A: \quad \fbox{1} & \fbox{2} & \fbox{3} & \fbox{4} \\
   \quad  \text{deck }B: \quad \fbox{1} & \fbox{4} & \fbox{2} & \fbox{3} \\
 \end{matrix}$$
 Once one of the above two transpositions occurs, card $\fbox{4}$ joins $\fbox{1}$ and $\fbox{3}$ as one of the cards that does label-to-label jumps only which will be unnoticed for the rest of the time in $[0,t_1]$.

 The other two options for card $\fbox{4}$ are label-to-label transpositions $<<\fbox{4},\fbox{1}>>$ and $<<\fbox{4},\fbox{2}>>$ getting
   $$\begin{matrix}
    \quad \text{deck }A: \quad \fbox{4} & \fbox{2} & \fbox{3} & \fbox{1} \\
   \quad  \text{deck }B: \quad \fbox{4} & \fbox{3} & \fbox{2} & \fbox{1} \\
 \end{matrix} ~~~\text{ and }~~~
  \begin{matrix}
    \quad \text{deck }A: \quad \fbox{1} & \fbox{4} & \fbox{3} & \fbox{2} \\
   \quad  \text{deck }B: \quad \fbox{1} & \fbox{3} & \fbox{4} & \fbox{2} \\
 \end{matrix}$$
 respectively, which we suppress noting.  The rates being the same for each one of $<\fbox{4},\fbox{3}>_A$,
 $<\fbox{4},\fbox{3}>_B$,   $<<\fbox{4},\fbox{1}>>$ and $<<\fbox{4},\fbox{2}>>$.
\item Card $\fbox{2}$ does the label to location jump $<<\fbox{2},i_1>>$ at time $t_1$.
\end{itemize}

\noindent
The point is that in Case II, either of the following transposition sequences
$$<\fbox{4},\fbox{3}>_A~~\text{ followed by }~~<<\fbox{2},i_1>>$$
or
$$<\fbox{4},\fbox{3}>_B~~\text{ followed by }~~<<\fbox{2},i_1>>$$
would lead to the discrepancies' cancelation.

If at time $t_1$ the discrepancies are not canceled out, start anew with a new random $i_1$ and exponential $t_1$.
There are two possibilities at these trials: either we end up with no discrepancies or we end up in the same boat as before and we try again.  Here we set the coupling rules by conditioning only on one event, $<<\fbox{2},i_1>>$. If we condition on more than one upcoming events (say  $<\fbox{4},\fbox{3}>_A$ in Case II), by conditioning inside the conditioning, we can increase the probability of
coincidence, thus producing a faster coupling time. Later in the paper we will deal with conditioning on a chain of events.
\vskip 0.2 in
\noindent
 The above example is a simplified version of the coupling construction to follow.

\subsubsection{Two discrepancies and one association map ($d=2,k=1$)}

   We will start with the case of two discrepancies, $d_1$ and $d_2$ as illustrated below, in (\ref{gluing0}), when the coupling speed must be the slowest, and introduce the first improvement to the classical coupling construction.
\begin{eqnarray} \label{gluing0}
   \begin{matrix}
    A_t: \quad \dots  & \fbox{4} & \fbox{6} & \fbox{\bf b} & \fbox{9} & \fbox{\bf a} & \fbox{8} & \fbox{\bf a1} & \fbox{2} & \dots   \\
    B_t: \quad \dots  & \fbox{4} & \fbox{6} & \fbox{\bf a} & \fbox{9} & \fbox{\bf b} & \fbox{8} & \fbox{\bf a1} & \fbox{2} & \dots \\
 \phantom{A_t: \quad \dots  }&  &  & \uparrow &  & \uparrow &  & \uparrow  &  & \\
 \phantom{A_t: \quad \dots  }&  &  & d_2 &  & d_1 &  & i_1 &  &
 \end{matrix}
\end{eqnarray}
 If we let the uncoupled cards $\fbox{\bf a}$ and $\fbox{\bf b}$ do the label-to-location transpositions, and the coupled cards do the label-to-label transpositions, where the transpositions are synchronized for the top and the bottom decks, $A_t$ and $B_t$. Then, because of the group invariance property that was mentioned earlier in the paper, the above configuration will not change until either of the uncoupled cards $\fbox{\bf a}$ or $\fbox{\bf b}$ transposes with one of the discrepancy locations $d_1$ or $d_2$, simultaneously in both the top and the bottom processes. The label-to-label transpositions of the uncoupled cards can be suppressed citing group invariance. If we do not adjust this coupling, the waiting time to cancel the two discrepancies will average ${n^2 \over 4}$, which is too large. Next, we introduce the first modification to the coupling construction.

 \vskip 0.3 in
 \noindent
  We can adjust the coupling as follows. We take one of the two uncoupled cards, for example $\fbox{\bf a}$. Randomly select a site $i_1$, and a random exponential time $t_1$, and condition on transposition ${\ll \fbox{\bf a},i_1 \gg}$ happening at time $t_1$. All with respect to group invariance. If $i_1=d_1$ or $d_2$, the discrepancies disappear, and the coupling is completed by letting all coupled cards do label-to-lable transpositions that we suppress to notice, citing the group invariance. We want to utilize ${\ll \fbox{\bf a},i_1 \gg}$ if $i_1\not= d_1$ or $d_2$.

 So we condition on the label-to-location transposition ${\ll \fbox{\bf a},i_1 \gg}$ occurring at time $t_1$. If $i_1\not= d_1$ or $d_2$, the following construction will help us set the transition rates for the coupling in the time interval $[0,t_1]$. The first step is to associate locations in the top and the bottom decks as follows:
 \begin{description}
  \item[(a)] we pair the site $i_1$ in $A_t$ with the site $d_1$ in $B_t$, and call it $i_1/d_1$,
  \item[(b)] we pair the site $d_2$ in $A_t$ with the site $i_1$ in $B_t$, and call it $d_2/i_1$
  \item[(c)] we pair the site $d_1$ in $A_t$ with the site $d_2$ in $B_t$, and call it $d_1/d_2$
\end{description}
Therefore, at time $t=0$, the decks are aligned accordingly:
\begin{eqnarray} \label{gluing1}
     \begin{matrix}
    A_t: \quad \dots  & \fbox{4} & \fbox{6} & \fbox{\bf b} & \fbox{9} & \fbox{\bf a1} & \fbox{8} & \fbox{\bf a} & \fbox{2} & \dots   \\
    B_t: \quad \dots  & \fbox{4} & \fbox{6} & \fbox{\bf a1} & \fbox{9} & \fbox{\bf b} & \fbox{8} & \fbox{\bf a} & \fbox{2} & \dots \\
 \phantom{A_t: \quad \dots  }&  &  & \uparrow &  & \uparrow &  & \uparrow  &  & \\
 \phantom{A_t: \quad \dots  }&  &  & d_2/i_1 &  & i_1/d_1 &  & d_1/d_2 &  &
 \end{matrix}
 \end{eqnarray}
 The above diagram reads as follows. The card $\fbox{\bf b}$ in the upper deck is located at the site $d_2$, while the card $\fbox{\bf a1}$ in the lower deck is located at the site $i_1$. Similarly, the card $\fbox{\bf a1}$ in the upper deck is located at the site $i_1$, while the card $\fbox{\bf b}$ in the lower deck is located at the site $d_1$. Last, the card $\fbox{\bf a}$ in the upper deck is located at the site $d_1$, while the card $\fbox{\bf a}$ in the lower deck is located at the site $d_2$. So the above configuration is exactly the same as in (\ref{gluing0}).

 \noindent
 What is different is that the transposition $\ll \fbox{\bf a},i_1 \gg$ in (\ref{gluing0}) is equivalent to the relabeling of the sites in (\ref{gluing1}) described in the following diagram:
 \vskip 0.2in
 \qquad \qquad \qquad \qquad\qquad \qquad \qquad \qquad
 \fbox{$\begin{matrix} d_1/d_2 \longrightarrow i_1\\ i_1/d_1 \longrightarrow d_1\\ d_2/i_1 \longrightarrow d_2\end{matrix}$}
 \vskip 0.2in
 \noindent
 The site association (\ref{gluing1}) will be called the association map with respect to transposition $\ll \fbox{\bf a},i_1 \gg$ and jump time $t_1$. From our new perspective, time $t_1$ is the time when the location names change according to the above rule.  We will say that the association map expires at $t_1$.

 \vskip 0.3 in
 \noindent
We use the above association map to set the rates for the new coupling process in the time interval $[0, t_1]$:

 \begin{itemize}
 \item  Rates for $\fbox{\bf a}$:  the first transposition $\ll \fbox{\bf a},i_1 \gg$ occurs at time $t_1$.
\vskip 0.2 in
 \item Rates for $\fbox{\bf a1}$:  $\fbox{\bf a1}$ does the label-to-location transpositions, where locations are  defined by the association
 map (\ref{gluing1}). In other words, transpositions $\ll \fbox{\bf a1},i \gg$ $~(i \not=i_1, d_1, d_2)~$, $\ll \fbox{\bf a1}, i_1/d_1 \gg$, $\ll \fbox{\bf a1}, d_2/i_1 \gg$ and $\ll \fbox{\bf a1}, d_1/d_2 \gg$ occur with the usual rate of ${1 \over n^2}$.

 \noindent
 Transposition $\ll \fbox{\bf a1}, i_1/d_1 \gg$ should be interpreted as simultaneous occurrence of transpositions $<\fbox{\bf a1}, i_1>_A$ and $<\fbox{\bf a1},d_1>_B$ in the top and the bottom decks. Similarly $\ll \fbox{\bf a1}, d_2/i_1 \gg$ should be interpreted as simultaneous occurrence of transpositions $<\fbox{\bf a1}, d_2>_A$ and $<\fbox{\bf a1},i_1>_B$ in the top and the bottom decks. Transposition $\ll \fbox{\bf a1}, d_1/d_2 \gg =\ll \fbox{\bf a1}, \fbox{\bf a} \gg$ is label-to-label and therefore should be suppressed.
\vskip 0.2 in
 \item Rates for $\fbox{\bf b}$:  transpositions $\ll \fbox{\bf b},d_1 \gg$ and $\ll \fbox{\bf b},d_2 \gg$, $< \fbox{\bf b},\fbox{\bf a1}>_A$, $< \fbox{\bf b},\fbox{\bf a1}>_B$ and $\ll \fbox{\bf b},i \gg$ $~(i \not=i_1, d_1, d_2)~$ occur independently with rate ${1 \over n^2}$ each.

 \item The rest of the cards do label-to-label jumps simultaneously in the top and the bottom process, which we suppress according to group invariance.
 \end{itemize}

 \noindent
We need to add one more condition. We observe that any of the following six transpositions combined with the occurrence of transposition $\ll \fbox{\bf a},i_1 \gg$ at time $t_1$ will lead to the discrepancies cancelation. The six transpositions are $\ll \fbox{\bf a1}, i_1/d_1 \gg$, $\ll \fbox{\bf a1}, d_2/i_1 \gg$, $\ll \fbox{\bf b},d_1 \gg$, $\ll \fbox{\bf b},d_1 \gg$, $< \fbox{\bf b},\fbox{\bf a1}>_A$,
 and $< \fbox{\bf b},\fbox{\bf a1}>_B$.

 \noindent
 Therefore, if one of the above mentioned six transpositions occurs before time $t_1$, from that time until time $t_1$, all the cards other than $\fbox{\bf a}$ will only be allowed to do label-to-label jumps.
\vskip 0.2 in

 \noindent
 Now, since we classified (and therefore suppressed) $\ll \fbox{\bf a1}, d_1/d_2 \gg$ as a label-to-label transposition, the jump time $t_2$ for the card $\fbox{\bf a1}$
 is exponential with rate ${n-1 \over n^2}$.
 Now, $P[t_2<t_1]={{n-1 \over n} \over 1+ {n-1 \over n}}={n-1 \over 2n-1}$, and
 conditioning on $t_2<t_1$, either $\ll \fbox{\bf a1}, i_1/d_1 \gg$ or $\ll \fbox{\bf a1}, d_2/i_1 \gg$ will occur with probability ${2 \over n-1}$. Once again observe that both will result in the discrepancy cancelation at time $t_1$.

 \noindent
  Let $t_b$ denote the first time one of the four transpositions $\ll \fbox{\bf b},d_1 \gg$, $\ll \fbox{\bf b},d_1 \gg$, $< \fbox{\bf b},\fbox{\bf a1}>_A$ and $< \fbox{\bf b},\fbox{\bf a1}>_B$ occurs. Then $t_b$ is exponential with rate ${4 \over n^2}$. Recall that if $t_b<t_1$ then the discrepancies will cancel out at time $t_1$.

 \noindent
 Therefore the probability of canceling the discrepancy by time $t_1$ is
 $$P[i_1=d_1~\text{or}~d_2]+\left(1-{2 \over n}\right)\cdot P[t_2<t_1]\cdot {2 \over n-1} + \left(1-{2 \over n}\left(2-{1 \over n-1}\right)\right)P[t_b<t_1] \approx {8 \over n},$$
 where $P[i_1=d_1~\text{or}~d_2]+\left(1-{2 \over n}\right)\cdot P[t_2<t_1]\cdot {2 \over n-1}={2 \over n}+\left(1-{2 \over n}\right)\cdot {n-1 \over 2n-1} \cdot {2 \over n-1}={2 \over n}\left(2-{1 \over n-1}\right)$ is the
 probability that the discrepancy is canceled by  $\ll \fbox{\bf a},i_1 \gg$ when $i=i_1$, or as the result of one of the two transpositions, $\ll \fbox{\bf a1}, i_1/d_1 \gg$ or $\ll \fbox{\bf a1}, d_2/i_1 \gg$, combined with $\ll \fbox{\bf a},i_1 \gg$. Now, if this does not happen, then with probability $P[t_b<t_1]={{4 \over n^2} \over {1 \over n}+{4 \over n^2}}\approx {4 \over n}$,
 the discrepancies can be canceled out with one of the four  transpositions, $\ll \fbox{\bf b},d_1 \gg$, $\ll \fbox{\bf b},d_1 \gg$, $< \fbox{\bf b},\fbox{\bf a1}>_A$
 or $< \fbox{\bf b},\fbox{\bf a1}>_B$, combined with $\ll \fbox{\bf a},i_1 \gg$.

 \noindent
 So, with probability $\approx {8 \over n}$, the discrepancies will cancel by time $t_1$. If they do not cancel, we repeat the association trick,
 thus coupling the two discrepancies in approximately ${n^2 \over 8}$ steps on average, instead of ${n^2 \over 4}$. However we can do better by using more than one association map at a time.

 \vskip 0.3 in
\noindent
{\bf Example.} Next, we illustrate how a cancelation of discrepancies on an association map (\ref{gluing1}) implies discrepancies' cancelation at time $t_1$ on (\ref{gluing0}). Consider the following case. If $t_2 < t_1$, and if $\ll \fbox{\bf a1}, d_2/i_1 \gg$ occurs  at time $t_2$, then we will observe the following dynamics on the association map.

\noindent
 The configurations will evolve from
  $$\begin{matrix}
    A_t: \quad \dots  & \fbox{4} & \fbox{6} & \fbox{\bf b} & \fbox{9} & \fbox{\bf a1} & \fbox{8} & \fbox{\bf a} & \fbox{2} & \dots   \\
    B_t: \quad \dots  & \fbox{4} & \fbox{6} & \fbox{\bf a1} & \fbox{9} & \fbox{\bf b} & \fbox{8} & \fbox{\bf a} & \fbox{2} & \dots \\
 \phantom{A_t: \quad \dots  }&  &  & \uparrow &  & \uparrow &  & \uparrow  &  & \\
 \phantom{A_t: \quad \dots  }&  &  & d_2/i_1 &  & i_1/d_1 &  & d_1/d_2 &  &
 \end{matrix}$$
 to
  $$\begin{matrix}
    A_t: \quad \dots  & \fbox{4} & \fbox{6} & \fbox{\bf a1} & \fbox{9} & \fbox{\bf b} & \fbox{8} & \fbox{\bf a} & \fbox{2} & \dots   \\
    B_t: \quad \dots  & \fbox{4} & \fbox{6} & \fbox{\bf a1} & \fbox{9} & \fbox{\bf b} & \fbox{8} & \fbox{\bf a} & \fbox{2} & \dots \\
 \phantom{A_t: \quad \dots  }&  &  & \uparrow &  & \uparrow &  & \uparrow  &  & \\
 \phantom{A_t: \quad \dots  }&  &  & d_2/i_1 &  & i_1/d_1 &  & d_1/d_2 &  &
 \end{matrix}$$
  at time $t_2$, and
   $$\begin{matrix}
    A_t: \quad \dots  & \fbox{4} & \fbox{6} & \fbox{\bf a1} & \fbox{9} & \fbox{\bf b} & \fbox{8} & \fbox{\bf a} & \fbox{2} & \dots   \\
    B_t: \quad \dots  & \fbox{4} & \fbox{6} & \fbox{\bf a1} & \fbox{9} & \fbox{\bf b} & \fbox{8} & \fbox{\bf a} & \fbox{2} & \dots \\
 \phantom{A_t: \quad \dots  }&  &  & \uparrow &  & \uparrow &  & \uparrow  &  & \\
 \phantom{A_t: \quad \dots  }&  &  & d_2 &  & d_1 &  & i_1 &  &
 \end{matrix}$$
 at time $t_1$.

\vskip 0.2 in
\noindent
 The above were the transformations one would see on the association map.  The corresponding
 evolutions of the decks with respect to the original site associations (\ref{gluing0}) will be as follows.
 From
 $$\begin{matrix}
    A_t: \quad \dots  & \fbox{4} & \fbox{6} & \fbox{\bf b} & \fbox{9} & \fbox{\bf a} & \fbox{8} & \fbox{\bf a1} & \fbox{2} & \dots   \\
    B_t: \quad \dots  & \fbox{4} & \fbox{6} & \fbox{\bf a} & \fbox{9} & \fbox{\bf b} & \fbox{8} & \fbox{\bf a1} & \fbox{2} & \dots \\
 \phantom{A_t: \quad \dots  }&  &  & \uparrow &  & \uparrow &  & \uparrow  &  & \\
 \phantom{A_t: \quad \dots  }&  &  & d_2 &  & d_1 &  & i_1 &  &
 \end{matrix}$$
 to
   $$\begin{matrix}
    A_t: \quad \dots  & \fbox{4} & \fbox{6} & \fbox{\bf a1} & \fbox{9} & \fbox{\bf a} & \fbox{8} & \fbox{\bf b} & \fbox{2} & \dots   \\
    B_t: \quad \dots  & \fbox{4} & \fbox{6} & \fbox{\bf a} & \fbox{9} & \fbox{\bf b} & \fbox{8} & \fbox{\bf a1} & \fbox{2} & \dots \\
 \phantom{A_t: \quad \dots  }&  &  & \uparrow &  & \uparrow &  & \uparrow  &  & \\
 \phantom{A_t: \quad \dots  }&  &  & d_2 &  & d_1 &  & i_1 &  &
 \end{matrix}$$
 at time $t_2$, and to
   $$\begin{matrix}
    A_t: \quad \dots  & \fbox{4} & \fbox{6} & \fbox{\bf a1} & \fbox{9} & \fbox{\bf b} & \fbox{8} & \fbox{\bf a} & \fbox{2} & \dots   \\
    B_t: \quad \dots  & \fbox{4} & \fbox{6} & \fbox{\bf a1} & \fbox{9} & \fbox{\bf b} & \fbox{8} & \fbox{\bf a} & \fbox{2} & \dots \\
 \phantom{A_t: \quad \dots  }&  &  & \uparrow &  & \uparrow &  & \uparrow  &  & \\
 \phantom{A_t: \quad \dots  }&  &  & d_2 &  & d_1 &  & i_1 &  &
 \end{matrix}$$
 at time $t_1$.

\subsubsection{Two discrepancies and $\varepsilon n$ association maps ($d=2,k=\lfloor \varepsilon n \rfloor $)}
Next we develop a calculus of association maps to show the coupling construction for the case of two discrepancies with coupling time of order $O(n \log{n})$.

\vskip 0.2 in
\noindent
Since the rate of ${n^2 \over 8}$ is still not good enough, we need to revise and enhance the coupling construction by introducing chains of association maps.
The first part is similar to the preceding subsection. Namely, we begin by randomly selecting a site $i_1$, and a random exponential time $t_1$, and conditioning on transposition ${\ll \fbox{\bf a},i_1 \gg}$ happening at time $t_1$. See (\ref{gluing0}).
\vskip 0.2 in
\noindent
In order to set the rates for $\fbox{\bf a1}$ in the coupling process, we considered the following  association map
  \begin{eqnarray} \label{gluing2}
 \begin{matrix}
    A_t: \quad \dots  & \fbox{\bf a2} & \fbox{6} & \fbox{\bf b} & \fbox{9} & \fbox{\bf a1} & \fbox{8} & \fbox{\bf a} & \fbox{2} & \dots   \\
    B_t: \quad \dots  & \fbox{\bf a2} & \fbox{6} & \fbox{\bf a1} & \fbox{9} & \fbox{\bf b} & \fbox{8} & \fbox{\bf a} & \fbox{2} & \dots \\
 \phantom{A_t: \quad \dots  }& \uparrow &  & \uparrow &  & \uparrow &  & \uparrow  &  & \\
 \phantom{A_t: \quad \dots  }& i_2 &  & d_2^* &  & d_1^* &  & i_1^* &  &
 \end{matrix}
 \end{eqnarray}
 where
 \begin{description}
  \item[(a)] $d^*_1$ denotes $i_1/d_1$ before $t_1$, and $d_1$ after $t_1$,
  \item[(b)] $d^*_2$ denotes $d_2/i_1$ before $t_1$, and $d_2$ after $t_1$,
  \item[(c)] $i^*_1$ denotes $d_1/d_2$ before $t_1$, and $i_1$ after $t_1$.
\end{description}
 On the above association map, we again have two discrepancies, and we adjust the coupling rules so that to cancel these new discrepancies. We do this by selecting a random  location $i_2 \in \{1,2,\dots, d_1^*,\dots,d_2^*, \dots,n\} \backslash \{i_1^*\}$ on the association map and a random time $t_2$, distributed exponentially ${1 \over n}\cdot \left(1-{1 \over n}\right)$. We condition on
the  transposition $\ll \fbox{\bf a1}, i_2 \gg$ happening at the time $t_2$. The group invariance allows us to suppress $\ll \fbox{\bf a1}, i_1^* \gg$ as   $\ll \fbox{\bf a1}, i_1^* \gg=\ll \fbox{\bf a1}, \fbox{\bf a} \gg$ is label-to-label.
\vskip 0.2 in
\noindent
 On the new scheme, if $i_2=d^*_1$ or $d^*_2$, the discrepancies cancel out, and the process couples after time $\tau_2=\max\{t_1,t_2\}$. See the example in the preceding subsection.
 \vskip 0.2 in
\noindent
 Now, if $i_2 \not= d^*_1$ or $d^*_2$, we will construct one more association map in order to set the transition rates for the card $\fbox{\bf a2}$ located at $i_2$. We pair the site $i_2$ in $A_t$ with site $d^*_1$ in $B_t$ and call it $i_2/d^*_1$; pair the site $d^*_2$ in $A_t$ with site $i_2$ in $B_t$ and call it $d^*_2/i_2$; and finally, pair the site $d^*_1$ in $A_t$ with site $d^*_2$ in $B_t$ and call it $d^*_1/d^*_2$.  After pairing anew the locations on top with the locations at the bottom, the association map (\ref{gluing2}) is realigned as
 \begin{eqnarray*}
   \begin{matrix}
    A_t: \quad \dots  & \fbox{\bf a1} & \fbox{6} & \fbox{\bf b} & \fbox{9} & \fbox{\bf a2} & \fbox{8} & \fbox{\bf a} & \fbox{2} & \dots   \\
    B_t: \quad \dots  & \fbox{\bf a1} & \fbox{6} & \fbox{\bf a2} & \fbox{9} & \fbox{\bf b} & \fbox{8} & \fbox{\bf a} & \fbox{2} & \dots \\
 \phantom{A_t: \quad \dots  }& \uparrow &  & \uparrow &  & \uparrow &  & \uparrow  &  & \\
 \phantom{A_t: \quad \dots  }& d^*_1/d^*_2 &  & d^*_2/i_2 &  & i_2/d^*_1 &  & i^*_1 &  &
 \end{matrix}
 \end{eqnarray*}
In other word, we construct the second association map
 \begin{eqnarray} \label{gluing3}
   \begin{matrix}
    A_t: \quad \dots  & \fbox{\bf a1} & \fbox{6} & \fbox{\bf b} & \fbox{9} & \fbox{\bf a2} & \fbox{8} & \fbox{\bf a}  & \fbox{\bf a3} & \dots   \\
    B_t: \quad \dots  & \fbox{\bf a1} & \fbox{6} & \fbox{\bf a2} & \fbox{9} & \fbox{\bf b} & \fbox{8} & \fbox{\bf a}  & \fbox{\bf a3} & \dots \\
 \phantom{A_t: \quad \dots  }& \uparrow &  & \uparrow &  & \uparrow &  & \uparrow  & \uparrow  & \\
 \phantom{A_t: \quad \dots  }& i^{**}_2 &  & d^{**}_2 &  & d^{**}_1 &  & i^*_1 & i_3  &
 \end{matrix}
 \end{eqnarray}
 where 
  \begin{description}
  \item[(a)] $d^{**}_1$ denotes $i_2/d^*_1$ before $t_2$, and $d^*_1$ after $t_2$,
  \item[(b)] $d^{**}_2$ denotes $d^*_2/i_2$ before $t_2$, and $d^*_2$ after $t_2$,
  \item[(c)] $i^{**}_2$ denotes $d^*_1/d^*_2$ before $t_2$, and $i_2$ after $t_2$.
\end{description}
 \vskip 0.2 in
\noindent
 Next, we pick a random location $i_3 \in \{1,2,\dots,d^{**}_1,\dots,d^{**}_2,\dots, n\} \backslash \{ i^*_1, i^{**}_2\}$ on the association map (\ref{gluing3}).
 We pick a random waiting time $t_3$, distributed exponentially with the parameter ${1 \over n}\cdot \left(1-{2 \over n}\right)$,  and condition on  transposition $\ll \fbox{\bf a2}, i_3 \gg$ on the association map (\ref{gluing3}) at time $t_3$. 
 Once again the group invariance allows us to suppress $\ll \fbox{\bf a2}, i_1^* \gg$ and $\ll \fbox{\bf a2}, i_2^{**} \gg$ as   $\ll \fbox{\bf a2}, i_1^* \gg=\ll \fbox{\bf a2}, \fbox{\bf a} \gg$ and $\ll \fbox{\bf a2}, i_2^{**} \gg=\ll \fbox{\bf a2}, \fbox{\bf a1} \gg$  are label-to-label.
 
 If $i_3=d^{**}_1$ or $d^{**}_2$, the discrepancies cancel out on the association map (\ref{gluing3}), and the process couples after time $\tau_3=\max\{t_1,t_2,t_3\}$. If $i_3 \not= d^{**}_1$ or $d^{**}_2$, we will construct the third association map in order to set up the transition rates for the card $\fbox{\bf a3}$ located at $i_3$. Proceeding inductively, after $j$ iterations,  we either cancel the discrepancies with probability of at least ${2 \over n-j+1}$, or construct one more association map.
 \noindent
 \vskip 0.2 in
 We construct a chain of at most $k=\lfloor \varepsilon n \rfloor $ association maps like that, where $\varepsilon \in (0,1)$ is fixed, thus setting the transition rates for the cards  $\fbox{\bf a}$, $\fbox{\bf a1}$, $\fbox{\bf a2},\dots,\fbox{\bf ak}$.
 As for all the cards different (w.r.t. group invariance) from $\fbox{\bf a}$, $\fbox{\bf a1}$, $\fbox{\bf a2},\dots,\fbox{\bf ak}$ and $\fbox{\bf b}$, all their jumps are set as label-to-label transpositions, which we suppress. 
  \vskip 0.2 in
  \noindent
Let $T_a$ be the first time one of the transpositions 
$$\ll \fbox{\bf a}, d_1 \gg ~\text{ and }~\ll \fbox{\bf a}, d_2 \gg,~~~\text{ and }~\ll \fbox{\bf aj}, d_1^{j*} \gg~\text{ and }~\ll \fbox{\bf aj}, d_2^{j*} \gg~\text{ for all }j=1,\dots,k$$
occurs. Then $T_a$ is exponential with rate ${2(k+1) \over n^2}$. Observe that each of the above transpositions leads to the discrepancies cancelation on of the $k$ association maps, or in the original association of sites.
   \vskip 0.2 in
  \noindent
 Lastly,  we have to set the rates for $\fbox{\bf b}$. We let the card $\fbox{\bf b}$ make label-to-label jumps
 $\ll \fbox{\bf b}, \fbox{\bf c} \gg$ simultaneously in the top in the bottom decks with rate ${1 \over n^2}$ whenever the card $\fbox{\bf c}$  is {\bf not} $\fbox{\bf a}$, or $\fbox{\bf a1}$, or $\fbox{\bf a2},\dots,\fbox{\bf ak}$. 
 As for transposing card $\fbox{\bf b}$ with $\fbox{\bf a}$, $\fbox{\bf a1}$, $\fbox{\bf a2},\dots,\fbox{\bf ak}$, we let the transpositions
  $$<\fbox{\bf b},\fbox{\bf a}>_A,<\fbox{\bf b},\fbox{\bf a}>_B, <\fbox{\bf b},\fbox{\bf a1}>_A,<\fbox{\bf b},\fbox{\bf a1}>_B,\dots,
  <\fbox{\bf b},\fbox{\bf ak}>_A,<\fbox{\bf b},\fbox{\bf ak}>_B$$
  occur independently in either the top or the bottom decks with the standard rate of ${1 \over n^2}$, up until the two discrepancies get canceled. 
  \vskip 0.2 in
  \noindent  
  If $T_b$ is the first time one of the label-to-label jumps 
  $$<\fbox{\bf b},\fbox{\bf a}>_A,<\fbox{\bf b},\fbox{\bf a}>_B, <\fbox{\bf b},\fbox{\bf a1}>_A,<\fbox{\bf b},\fbox{\bf a1}>_B,\dots,
  <\fbox{\bf b},\fbox{\bf ak}>_A,<\fbox{\bf b},\fbox{\bf ak}>_B$$
 occurs, then $T_b$ must be exponential with rate ${2(k+1) \over n^2}$. Each one of these jumps leads to the discrepancies cancelation on one of the $k$ association maps, or in the original association of sites.
Thus the average time for the discrepancies cancelation on at least one of the association maps will be bounded above by
  $$E[T_a \wedge T_b] ={n^2  \over 4(k+1)} \leq {n  \over 4\varepsilon}$$

   \noindent
  Now, the discrepancies cancelation on anyone of the $k$ association map leads to cancelation for the original association of locations, and coupling by the time $\tau_k=\max\{t_1,t_2,\dots,t_k\}$ when each association map expires. Here for each $j \in \{1,2,\dots,k\}$,
 \vskip 0.2 in
  \qquad  \qquad  \qquad  \qquad  \qquad {$t_j$ is exponential with rate ${1 \over n}\cdot \left(1-{j-1 \over n}\right)$},
  \vskip 0.2 in
 \noindent
 and for $n$ large enough,
 $$E[\tau_k] \leq {1 \over {1 \over n}\cdot [1-{k-1 \over n}]} \log{k}  \leq {n \over 1-\varepsilon} \log{(\varepsilon  n)}.$$
 The $k$ association maps will, on average, expire in less than ${n \over 1-\varepsilon} \log{(\varepsilon n)}$ units of time.
  \vskip 0.2 in
  \noindent  
Thus the upper bound of ${n \over 1-\varepsilon} \log{(\varepsilon n)}+{n  \over 4\varepsilon}$ on the expected coupling time $E[T_c]$ is obtained in the case of two discrepancies.

\subsection{The general case of $d\geq 2$ discrepancies}
In this subsection we consider the situation when there are  $d \geq 2$ discrepancies.  For each of the $d$ discrepancy cards \fbox{${\bf b_{m}}$} ($m=1,\dots,d$), we condition on the first label-to-location jump $\ll\fbox{\bf b\phantom{g}}\!\!\!\!\!_{\bf m}~,i_{m,1} \gg$, and use this to set the transition rates for the cards \fbox{${\bf a_{m,1}}$} occupying the locations  $i_{m1}$ at time zero, all with respect to the group invariance. We do this by creating corresponding association maps. Next, we condition on the jumps  $\ll$\fbox{${\bf a_{m,1}}$},$i_{m,2} \gg$, where $i_{m,2}$ ($m=1,\dots,d$) are locations on the corresponding association maps. We use this in order to set the rates for the cards \fbox{${\bf a_{m,2}}$} occupying the locations  $i_{m,2}$ at time zero. Thus for each $m$ we create a chain of conditionings  
\fbox{${\bf b_m}$} $\rightarrow$ \fbox{${\bf a_{m,1}}$} $\rightarrow$ \fbox{${\bf a_{m,2}}$}$\rightarrow \dots$$\rightarrow$ \fbox{${\bf a_{m,j}}$}$\rightarrow \dots$
We let the rest of the cards do label-to-label jumps, simultaneously in the top and the bottom decks.

  \vskip 0.2 in
  \noindent 
  Note on the notations: in this subsection,  we will not use the $*$ symbol in the notations of discrepancies on the association maps.
    \vskip 0.2 in
  \noindent 
  In the picture example below, $d_1, d_2 \dots, d_d$ denote all the discrepancies, and  \fbox{${\bf b_1}$}, \fbox{${\bf b_2}$}, \dots, \fbox{${\bf b_d}$} denote all the discrepancy cards.
  $$\begin{matrix}
    A_t: \quad \dots  & \fbox{4} & \fbox{6} & \fbox{\bf b\scriptsize{2}} & \fbox{9} & \fbox{\bf b\scriptsize{8}} & \fbox{\bf b\scriptsize{3}} & \fbox{\bf b\scriptsize{4}} & \fbox{\bf b\scriptsize{1}} & \dots   \\
    B_t: \quad \dots  & \fbox{4} & \fbox{6} & \fbox{\bf b\scriptsize{8}} & \fbox{9} & \fbox{\bf b\scriptsize{1}} & \fbox{\bf b\scriptsize{4}} & \fbox{\bf b\scriptsize{2}} & \fbox{\bf b\scriptsize{3}} & \dots \\
 \phantom{A_t: \quad \dots  }&  &  & \uparrow &  & \uparrow &  \uparrow &  \uparrow &  \uparrow & \\
 \phantom{A_t: \quad \dots  }&  &  & d_1 &  & d_2 & d_3 & d_4 & d_5 &
 \end{matrix}$$
  In the general case of $d\geq 2$ discrepancies, for each discrepancy card  \fbox{${\bf b_1}$}, \fbox{${\bf b_2}$}, \dots, \fbox{${\bf b_d}$} we construct a chain of association maps via conditioning, adding new association maps at expiration times.   Each chain will determine the rates for
  \fbox{${\bf b_m}$},  \fbox{${\bf a_{m,1}}$}, \fbox{${\bf a_{m,2}}$},...~.
  A transposition $\ll$\fbox{${\bf a_{m,j}}$},$d_i\gg$ will cancel discrepancies on the $j$th association map, in the $m$th chain of association  maps.
  For each $m$ from $1$ to $d$,  a transposition  $\ll\fbox{\bf b\phantom{g}}\!\!\!\!\!_{\bf m}~,i_{m1} \gg$
  is set to take place at time $t_{m,1}$, and conditioning on this, an association map is created.  Here again, $i_{2,1} \not=i_{1,1}$, $i_{3,1} \not=i_{1,1}, i_{2,1}$ etc., and therefore $t_{1,1}, t_{2,1}, \dots$  are exponential random variables with respective rates ${1 \over n}, {n-1 \over n^2}, {n-2 \over n^2}, \dots$.
  Transposition $\ll$\fbox{${\bf a_{m,1}}$},$i_{m,2} \gg$,   where $i_{m,2}$ ($i_{m,2} \not= i_{1,1}, \dots,i_{d,1}, i_{2,1},\dots,i_{(m-1),2}$) is a location with respect the corresponding association map, is conditioned to take place at time $t_{m2}$, and a new association  map is created. Since all locations $i_{mj}$ corresponding to the $k=\lfloor \varepsilon n \rfloor $ association maps are different,  for a given fixed $\kappa \in (\varepsilon, 1)$, if we keep the number of conditionings below $\kappa n$ at all times, then each expiration time $t_{mj}$ will be assigned a Poisson exponential rate of at most ${1-\kappa \over n}$. 

\vskip 0.2 in 
\noindent
Observe that in the $d=2$ case considered in the preceding subsection, the main point of the construction was that without the adjustment, only the discrepancy cards could cancel the discrepancies. There, using the conditionings and association maps we adjusted the coupling to enable discrepancy cancelations via the jumps by the non-discrepancy cards \fbox{${\bf a1}$}, \fbox{${\bf a2}$}, \dots, \fbox{${\bf ak}$} to discrepancy locations on the corresponding association maps. Here too, we want the non-discrepancy cards to cancel the discrepancies. We condition on the jumps of  \fbox{${\bf b_1}$}, \fbox{${\bf b_2}$}, \dots, \fbox{${\bf b_d}$} to make  \fbox{${\bf a_{1,1}}$}, \fbox{${\bf a_{2,1}}$}, \dots, \fbox{${\bf a_{d,1}}$} the new discrepancy cards on the corresponding association maps, enabling them to cancel the discrepancies by jumping to a discrepancy location. 
 $$\begin{matrix}
    A_t: \quad \dots  & \fbox{4} & \fbox{6} & \fbox{\bf a\tiny{2,1}} & \fbox{9} & \fbox{\bf a\tiny{8,1}} & \fbox{\bf a\tiny{3,1}} & \fbox{\bf a\tiny{4,1}} & \fbox{\bf a\tiny{1,1}} & \dots   \\
    B_t: \quad \dots  & \fbox{4} & \fbox{6} & \fbox{\bf a\tiny{8,1}} & \fbox{9} & \fbox{\bf a\tiny{1,1}} & \fbox{\bf a\tiny{4,1}} & \fbox{\bf a\tiny{2,1}} & \fbox{\bf a\tiny{3,1}} & \dots \\
 \phantom{A_t: \quad \dots  }&  &  & \uparrow &  & \uparrow &  \uparrow &  \uparrow &  \uparrow & \\
 \phantom{A_t: \quad \dots  }&  &  & d_1 &  & d_2 & d_3 & d_4 & d_5 &
 \end{matrix}$$
\noindent
We condition on the jumps of \fbox{${\bf a_{1,1}}$}, \fbox{${\bf a_{2,1}}$}, \dots, \fbox{${\bf a_{d,1}}$} to make \fbox{${\bf a_{1,2}}$}, \fbox{${\bf a_{2,2}}$}, \dots, \fbox{${\bf a_{d,2}}$} the discrepancy cards on the new association maps, and so on. Once a discrepancy is canceled out on an association map, the card that was coupled will not be used in the consequent rounds of conditionings.


 \vskip 0.2 in
 \noindent 
 Therefore, the average time for the discrepancy cancelation
  over all $k$ association maps is 
  $$E[T_a] \leq {n^2  \over kd} \leq {n  \over \varepsilon d}$$
  Since we are dealing with random transpositions, the construction can be adjusted to couple  twice as fast. So without loss of generality we expect the first discrepancy cancelation time to be  $\leq {n  \over 2\varepsilon d}$.
  \vskip 0.3 in

  \noindent

  When a discrepancy is canceled, after an average of ${n \over 2\varepsilon d}$ units of time, we will reuse all of the $d-1$ chains of association maps that were not responsible for the cancelation.  So at most, $\lceil {\varepsilon n \over d} \rceil$ association maps will not be reused.  Every time we cancel one of $d$ discrepancies, we are left with an average of
  ${1 \over 2}\cdot {\varepsilon n \over d}$ non-reusable association maps in the corresponding chain as chain of association maps  has length $\approx {\varepsilon n \over d}$. 

 If we impose an upper bound $\kappa n$ ($\varepsilon \leq \kappa <1$) on the number of association maps allowed at any one time,
  then the average time of discrepancy cancelation will be delayed by at most
  $${n \over 1-\kappa} \cdot \log{\left(1+{\varepsilon \over (\kappa -\varepsilon)d}\right)}
  \leq {\varepsilon \over (1-\kappa)(\kappa -\varepsilon)} \cdot {n \over d}~.$$
  Thus it will take less than
  $$ \left[ {1\over 2\varepsilon}+{\varepsilon \over (1-\kappa)(\kappa -\varepsilon)}\right] \cdot n \log{n}$$
  time steps to cancel every discrepancy, each w.r.t. some association map.

   After the discrepancies cancel with respect to the association maps, it will take an average of about ${n \over 1-\kappa} \log{(\kappa n)}$
  units of time for all remaining association maps to expire.
  Thus the upper bound on the expectation of coupling time will be
  $$\left[ {1\over 2\varepsilon}+{\kappa \over (1-\kappa)(\kappa -\varepsilon)}\right]  \cdot n \log{n}~~\text{ for any }~~ 0<\varepsilon < \kappa <1.$$
 Thus we obtained a $C n\log{n}$ upper bound on coupling time, where after taking optimal $\varepsilon$ and $\kappa$, constant $C$ becomes fairly small. Namely $C <6$.
 Observe that in the above constructed upper bound, the term  ${1\over 2\varepsilon} n \log{n}$ is responsible for the discrepancy cancelations on the association maps. From the information theory perspective this points at the ${1\over 2} n \log{n}$ cut-off for the mixing time.

\bibliographystyle{amsplain}


\end{document}